\documentclass[10pt,twoside,english]{amsart}
\advance\oddsidemargin by -1.0cm
\advance\evensidemargin by -1.0cm
\advance\topmargin by -1.0cm 
\usepackage{amssymb}
\usepackage{babel}
\usepackage{amstext}
\usepackage{amscd}   
\usepackage{epsfig} 
 \usepackage{rotating}
\theoremstyle{plain}
\newtheorem{cor}{Corollary}[section]
\newtheorem{lem}{Lemma}[section]

\theoremstyle{definition}

\newtheorem{NB}{Remark}[section]

\newcommand{\R}{\ensuremath{\mathbb{R}}}

\newcommand{\Ric}{\ensuremath{\mathrm{Ric}}}

\newtheorem{proposition}{Proposition}[section]

\newtheorem{theorem}[proposition]{Theorem}

\newcommand{\be}{\begin{equation}}
\newcommand{\ee}{\end{equation}}
\newcommand{\bea}{\begin{eqnarray}}
\newcommand{\eea}{\end{eqnarray}}
\newcommand{\bean}{\begin{eqnarray*}}
\newcommand{\eean}{\end{eqnarray*}}

\begin{document}

\title{Curvature dependent lower bounds for the first eigenvalue of the
Dirac operator}
\author{K.-D. Kirchberg}
\date{}
\maketitle \vspace{1cm}

\begin{abstract}
\noindent Using Weitzenb\"ock techniques on any compact Riemannian spin
manifold we derive inequalities that involve a real parameter and 
join the eigenvalues of the Dirac ope\-ra\-tor with curvature terms. The 
discussion of these inequalities yields vanishing theorems for the kernel
of the Dirac ope\-ra\-tor $D$ and lower bounds for the spectrum of $D^2$
if the curvature satisfies certain conditions. \\

\end{abstract}
\bigskip

\noindent\section{\bf Introduction}

In 1980 Th. Friedrich \cite{1} proved that, on any compact Riemannian
spin $n$-manifold $M$ of scalar curvature $S$ with $S_0 := \min \{ S(x) 
| \ x \in M \} >0$, every eigenvalue $\lambda$ of the Dirac operator $D$
satisfies the inequality
\begin{equation} \label{gl-01}
\lambda^2 \ge \frac{n}{4(n-1)} S_0 .
\end{equation}

In special geometric situations, better estimates are known (see \cite{5},
\cite{7}). For example, if $M$ is a spin K\"ahler manifold of complex
dimension $m$ and scalar curvature $S >0$, we have the inequalities
\begin{equation} \label{gl-02}
\lambda^2 \ge \left\{ \begin{array}{l}
\frac{m+1}{4m} S_0 \quad \quad \mbox{($m$ odd)}\\[1em]
\frac{m}{4(m-1)} S_0 \quad \mbox{($m$ even)} \ . 
		      \end{array} \right.
\end{equation}

The estimates (1), (2) are sharp in the sense that there are
manifolds for which the given lower bound itself is an eigenvalue of
$D^2$. But this kind of estimate by the scalar curvature only is not useful
if $S$ has zeros or attains negative values. Hence, the question arises if 
there exist lower bounds for the spectrum of $D^2$ that depend on additional
curvature terms. For certain manifolds whose curvature tensor or Weyl tensor, 
respectively, is divergence-free (co-closed and, hence, harmonic) such
lower bounds have been obtained recently (see \cite{2}, \cite{3}). 
In the case of a 
compact Riemannian spin $n$-manifold $M$ with divergence-free curvature
tensor $R (\delta R=0)$, scalar curvature $S=0$, and nowhere vanishing Ricci
tensor, for example, the estimate
\begin{equation} \label{gl-03}
\lambda^2 > \frac{1}{4} \cdot \frac{|\Ric|_0^2}{| \kappa_0| + | \Ric |_0 
\sqrt{\frac{n-1}{n}}}
\end{equation}

\thispagestyle{empty} 

is valid, where $| \Ric |_0 >0$ denotes the minimum of the length of the
Ricci tensor and $\kappa_0$ the smallest eigenvalue of $\Ric$ on $M$
(\cite{2}, Th. 2.2). Moreover, it has been proved that $\ker (D)$ is
trivial, i.e., there are no harmonic spinors if $M$ is compact with 
divergence-free curvature tensor and scalar curvature $S \le 0$ such that
the inequality
\begin{equation} \label{gl-04}
| \Ric |^2_0 > S \cdot \kappa_0
\end{equation}

holds (\cite{2}, Th. 2.2). We recall that $S$ is constant here, 
since the supposition $\delta R=0$ is equivalent to the symmetry property
\begin{equation} \label{gl-05}
(\nabla_X \Ric )Y= (\nabla_Y \Ric )X
\end{equation}

of the covariant derivative $\nabla \Ric$ of the Ricci tensor,
 which immediately
implies $dS=0$. A more general supposition than (\ref{gl-05}) is
\begin{equation} \label{gl-06}
(\nabla_X \Ric)Y - (\nabla_Y \Ric)X= \frac{1}{2 (n-1)} (X(S)Y-
Y(S)X) . 
\end{equation}

For dimension $n \ge 4$, (\ref{gl-06}) is equivalent to the condition that
the Weyl tensor $W$ is divergence-free $(\delta W=0)$ and, hence, harmonic
$(dW=0, \delta W=0)$. In the compact conformally non-flat case with $\delta
W=0$, the estimate
\begin{equation} \label{gl-07}
\lambda^2 \ge \frac{1}{8(n-1)} \left( (2n-1) S_0 + \sqrt{S_0^2 + 
\frac{n-1}{n} (\frac{4 \nu_0}{\mu})^2} \right)
\end{equation}

was proved for any eigenvalue $\lambda$ of the Dirac operator, where $\nu_0
\ge 0$ and $\mu >0$ are conformal invariants depending on $W$ only. For $S_0
>0$, (\ref{gl-07}) yields a better estimate than (\ref{gl-01}) if $\nu_0
>0$. For $S_0 \le 0$, the lower bound in (\ref{gl-07}) is positive if
$2 \nu_0 > n \mu | S_0 |$ (\cite{3}, Th. 3.1). In this paper we prove 
estimates similar to (3) and (\ref{gl-07}) which,
however, do not make use of
the suppositions (\ref{gl-05}) or (\ref{gl-06}), respectively. Moreover,
we obtain vanishing theorems for the space $\ker (D)$ of harmonic spinors
which are generalizations of those in \cite{2} and \cite{3}. Our results
are based on Weitzenb\"ock formulas for modified twistor operators, which
can partially be found in \cite{2}, \cite{3}, \cite{6} already.
However, what is new in this paper is the combination
of the various Weitzenb\"ock formulas for the modified twistor operators.
\\[1em]

\section{\bf Curvature endomorphisms of the spinor bundle}

Let $M$ be any Riemannian spin $n$-manifold with Riemannian metric $g$ and 
spinor bundle $\Sigma$. As usual, we denote by $\nabla$ the covariant 
derivative induced by $g$ on vector fields as well as on spinor fields
(Levi-Civita connection). For any vector fields $X, Y, Z$ and any spinor field
$\psi$, the Riemannian curvature tensor $R$ and the corresponding curvature
tensor $C$ of the spinor bundle are defined by
\begin{displaymath}
R(X,Y)Z:= \nabla^2_{X,Y} Z- \nabla^2_{Y,X} Z \quad , \quad C(X,Y) \psi :=
\nabla^2_{X,Y} \psi - \nabla^2_{Y,X} \psi , 
\end{displaymath}

where we use the notation
\begin{displaymath}
\nabla^2_{X,Y} := \nabla_X \circ \nabla_Y - \nabla_{\nabla_X Y}
\end{displaymath}

for the tensorial derivatives of second order. Given a local frame of vector 
fields $(X_1 , \ldots , X_n)$, we denote by $(X^1 , \ldots , X^n)$ the 
associated coframe defined by $X^k := g^{kl} X_l$, where $(g^{kl})$ is the 
inverse of the matrix $(g_{kl})$ with $g_{kl} := g(X_k , X_l)$. Thus,
for any orthonormal frame, we have $X^k = X_k \ (k=1 , \ldots , n)$. 
Then the Ricci tensor $\Ric$, the scalar curvature $S$, and the Dirac operator
$D$ are locally given by $\Ric (X) =R(X, X_k) X^k, S= \mathrm{tr} (\Ric) =
g(\Ric (X_k) , X^k)$ and $D \psi = X^k \cdot \nabla_{X_k} \psi$, 
respectively.\\

For the reader's convenience, we summarize some well-known, 
important identities:
\begin{equation} \label{gl-08}
C(X,Y)= \frac{1}{4} X_k \cdot R(X,Y)X^k , 
\end{equation}
\begin{equation} \label{gl-09}
X_k \cdot C(X^k ,X)= \frac{1}{2} \Ric (X) =C(X^k , X) \cdot X_k , 
\end{equation}
\begin{equation} \label{gl-10}
X_k \cdot \Ric (X^k)= -S = \Ric (X^k ) \cdot X_k , 
\end{equation}
\begin{equation} \label{gl-11}
X^k \cdot \nabla^2_{X_k , X} \psi = \nabla_X D \psi +  \frac{1}{2} \Ric (X)
\cdot \psi , 
\end{equation}
\begin{equation} \label{gl-12}
X^k \cdot \nabla^2_{X,X_k } \psi = \nabla_X D \psi . 
\end{equation}

The curvature endomorphism $C(X,Y)$ is anti-selfadjoint with respect to the 
Hermitian scalar product $\langle \cdot , \cdot \rangle$ on $\Sigma$, i.e., 
we have
\begin{equation} \label{gl-13}
C(X,Y)^* = - C(X,Y) . 
\end{equation}

Thus, the endomorphism $C^2 (X,Y):= C(Y,X_k) \circ C(X^k , X)$ has the 
property
\begin{equation} \label{gl-14}
C^2 (X,Y)^* = C^2 (Y,X)
\end{equation}

and, hence, the endomorphism $G: = C^2 (X_k , X^k)$ of $\Sigma$ is selfadjoint
and nonnegative
\begin{equation} \label{gl-15}
G^* =G \quad , \quad G \ge 0 . 
\end{equation}

Let $W$ denote the Weyl tensor of $M$ and consider the curvature 
endomorphisms $B(X,Y):=  \frac{1}{4} X_k \cdot W(X,Y) X^k$, $ B^2 (X,Y):=
B(Y, X_k) \circ B(X^k , X) , H:= B^2 (X_k , X^k)$. Then we have analogously:
\begin{equation} \label{gl-16}
B(X,Y)^* = - B(X,Y) \quad , \quad B^2 (X, Y)^* =B^2 (Y, X) , 
\end{equation}
\begin{equation} \label{gl-17} X_k \cdot B(X^k , X)= 0= B( X^k , X) \cdot 
X_k , 
\end{equation}
\begin{equation} \label{gl-18}
H^* = H \quad , \quad H \ge 0 . 
\end{equation}

The following lemma is proved by straightforward calculations.\\

\begin{lem} \label{lem-2-1}
The endomorphisms $G$ and $H$ are related by 
\begin{equation} \label{gl-19}
G=H+ \frac{1}{8} (|R|^2- |W|^2)= H + \frac{1}{2(n-2)} |\Ric - \frac{S}{n} 
|^2 + \frac{S^2}{4n(n-1)} . 
\end{equation}

Moreover, if $H = H_0 + H_2 + H_4$ is the decomposition of $H$ in the 
Clifford algebra into the components $H_0, H_2, H_4$ of degree $0,2$ and $4$,
respectively, then
\begin{equation} \label{gl-20}
H_0 = \frac{1}{8} | W|^2 \quad , \quad  H_2 =0 .
\end{equation}
\end{lem}

\vspace{0.5cm}

Using the notations $\delta R(X) := (\nabla_{X_k} R)(X,X^k) , \delta C(X):=
(\nabla_{X_k} C)(X, X^k)$ and \\
 $\delta W(X) := (\nabla_{X_k} W)(X, X^k) , 
\delta B(X):=(\nabla_{X_k} B)(X,X^k)$ we have the equations
\begin{equation} \label{gl-21}
\delta C(X)= \frac{1}{4} X_k  \cdot \delta R(X)X^k \quad , \quad \delta B
(X) = \frac{1}{4} X_k \cdot \delta W(X)X^k . 
\end{equation}

Moreover, it holds that
\begin{equation} \label{gl-22} 
\delta B(X)= \delta C(X) + \frac{1}{8(n-1)} (X \cdot dS - dS \cdot X) . 
\end{equation}

The second Biancchi identity implies
\begin{equation} \label{gl-23}
g(\delta R(X) Y,Z) = g(( \nabla_Y \Ric )Z -(\nabla_Z \Ric )Y,X) . 
\end{equation}

Inserting this into (\ref{gl-21}) we obtain
\begin{equation} \label{gl-24}
\delta C(X)= \frac{1}{4} (X^k \cdot (\nabla_{X_k} \Ric)X - (
\nabla_{X_k} \Ric )X \cdot X^k) . 
\end{equation}

Using (\ref{gl-21}) and (\ref{gl-24}) we find the identities
\begin{equation} \label{gl-25}
X_k \cdot \delta C(X^k)= \frac{1}{4} dS \quad , \quad X_k \cdot
\delta B(X^k) =0 . 
\end{equation}

For any vector field $X$, the endomorphisms $\delta C(X)$ and $\delta B(X)$
of $\Sigma$ are antiselfadjoint
\begin{equation} \label{gl-26}
\delta C(X)^* = - \delta C(X) \quad , \quad \delta B(X)^* = - \delta B(X).
\end{equation}

Thus, the endomorphisms $E:= - \delta C(X_k) \circ \delta C(X^k)$ and 
$F:= - \delta B(X_k) \circ \delta B(X^k)$ are selfadjoint and nonnegative
\begin{equation} \label{gl-27}
E^* = E  \quad , \quad F^* =F \quad , \quad E \ge 0\quad , \quad F \ge 0 .
\end{equation}

By (\ref{gl-22}) and (\ref{gl-25}), we obtain
\begin{equation} \label{gl-28}
E=F + \frac{1}{16(n-1)} |dS |^2 . 
\end{equation}

Moreover, by Proposition 3.1. in \cite{6}, it holds that
\begin{equation} \label{gl-29}
E= \frac{1}{4} | \nabla \Ric |^2 - \frac{1}{16} |dS|^2 + \frac{1}{8}
[ \nabla_{X_j} \Ric , \nabla_{X_k} \Ric] (X_l) \cdot X^j \cdot X^k \cdot
X^l , 
\end{equation}

where $[ \cdot , \cdot ]$ denotes the commutator of endomorphisms. Now we
introduce some numbers that occur in our following eigenvalue estimates.
Let $M$ be compact. We denote by $\nu_0$ the infimum of all eigenvalues of
$H$ on $M$. By definition, $\nu_0$ is a conformal invariant and we have the
inequality
\begin{equation} \label{gl-30}
\nu_0 | \psi |^2 \le \langle H \psi , \psi \rangle 
\end{equation}

for any $\psi \in \Gamma (\Sigma)$. By (\ref{gl-19}) we see that $\Ric$ and 
$\nu_0$ are obstructions against the existence of parallel spinors since
$\nabla \psi =0$ implies $C(X,Y) \cdot \psi =0$ for all vector fields
$X,Y$ and, hence, $G \psi =0$. The Schr\"odinger-Lichnerowicz formula
\begin{equation} \label{gl-31}
\nabla^* \nabla =D^2 - \frac{S}{4}
\end{equation}

shows that, in the compact case with vanishing scalar curvature, any
harmonic spinor $\psi$ $(D \psi =0)$ is parallel. Hence, $\ker (D)=0$ follows
if $M$ is compact and Ricci flat, but $\nu_0 >0$. In special situations, $\nu_0$
can easily be computed (\cite{3}, Section 3). Further, we consider the 
number
\begin{displaymath}
\mu := \sup \{ \| B(X,Y) \| \Big| x \in M, X, Y \in T_x M, 
g(X,Y)=0 , |X| =|Y| =1 \} , 
\end{displaymath}

where $\| \cdot \|$ denotes the operator norm. By definition, $\mu \ge 0$ is
a conformal invariant. By $\zeta$ we denote the corresponding supremum if 
$B$ is replaced by the spin curvature tensor $C$.\\

\begin{lem} \label{lem-2-2}
For any $\psi \in \Gamma (\Sigma)$, the inequalities
\begin{equation} \label{gl-32}
| \langle C^2 (X^k , X^l) \cdot \nabla_{X_k} \psi , \nabla_{X_l} 
\psi \rangle | \ \le (n-1)^2 \zeta^2 | \nabla \psi |^2 , 
\end{equation}
\begin{equation} \label{gl-33}
| \langle B^2 (X^k , X^l) \cdot \nabla_{X_k} \psi , \nabla_{X_l} \psi 
\rangle | \ \le (n-1)^2 \mu^2 | \nabla \psi |^2
\end{equation}

are valid.
\end{lem}
\vspace{0.3cm}

\begin{proof}
Let $(X_1 , \ldots , X_n)$ be any local orthonormal frame. 
Then, for all $k,l \in \{ 1 , \ldots , n \}$, we have the estimate\\

$\displaystyle (*) \quad \hfill \sum\limits^n_{j=1} \| C(X_j , X_k) \|
\ \| C(X_j , X_l) \| \le \left\{ \begin{array}{cl}
(n-1) \zeta^2 & \mbox{if $k=l$}\\
(n-2) \zeta^2 & \mbox{if $k \not= l$ . }
				 \end{array}
\right\} \hfill \mbox{}$\\

Now it holds that
\begin{eqnarray*}
&& | \langle C^2 (X^k , X^l) \nabla_{X_k} \psi , \nabla_{X_l} \psi \rangle
| \le \sum\limits_{j,k,l} | C(X_j , X_k) \nabla_{X_k} \psi , 
C(X_j , X_l) \nabla_{X_l} \psi \rangle |
\end{eqnarray*}
\begin{eqnarray*}
& \le & \sum\limits_{j,k,l} | C (X_j , X_k) \nabla_{X_k} \psi |  \ |
C(X_j , X_l) \nabla_{X_l} \psi | \le \sum\limits_{j,k,l} \| C(X_j , X_k) 
\| \ \| C(X_j , X_l)\| \ | \nabla_{X_k} \psi | | \nabla_{X_l} \psi | \\
&=& \sum\limits_{j,k} \| C(X_j , X_k) \|^2 | \nabla_{X_k} \psi |^2 +
\sum\limits_{j,k \not= l} \| C(X_j , X_k) \| \ \| C(X_j , X_l) \| \  | 
\nabla_{X_k} \psi |  | \nabla_{X_l} \psi | \\
& \stackrel{(*)}{\le} & (n-1) \zeta^2 | \nabla \psi |^2 + (n-2) \zeta^2
\sum\limits_{k \not= l} | \nabla_{X_k} \psi | | \nabla_{X_l} \psi | \\
&=& \zeta^2 | \nabla \psi |^2 + (n-2) \zeta^2 \sum\limits_{k,l} | 
\nabla_{X_k} \psi| | \nabla_{X_l} \psi |\\
&=& \zeta^2 | \nabla \psi |^2 + (n-2) \zeta^2 ( \sum\limits_k | 
\nabla_{X_k} \psi | )^2 \\
& \le & \zeta^2 | \nabla \psi |^2 + n(n-2) \zeta^2 \sum\limits_k | \nabla_{X_k}
\psi |^2 =( n-1)^2 \zeta^2 | \nabla \psi |^2 . 
\end{eqnarray*}

This proves (\ref{gl-32}). An analogous calculation yields (\ref{gl-33}).
\end{proof}

\vspace{0.3cm}

We remark that (\ref{gl-33}) is a better estimate than the corresponding 
estimate (\ref{gl-23}) in \cite{3}. \\[1em]

\section{\bf Estimates depending on the Ricci tensor}

\newcommand{\cald}{\mathcal{D}}
\newcommand{\calp}{\mathcal{P}}
\newcommand{\cale}{\mathcal{E}}
\newcommand{\calr}{\mathcal{R}}
\newcommand{\cals}{\mathcal{S}}
\newcommand{\calt}{\mathcal{T}}
\newcommand{\calq}{\mathcal{Q}}

Let $M$ be a Riemannian spin $n$-manifold and 
\begin{displaymath}
\cald : \Gamma (\Sigma) \to \Gamma (TM \otimes \Sigma)
\end{displaymath}

the corresponding twistor operator locally given by $\cald \psi := X^k
\otimes \cald_{X_k} \psi$ with 
\begin{displaymath}
\cald_X \psi := \nabla_X \psi + 
\frac{1}{n} X \cdot D \psi . 
\end{displaymath}

For $s,t \in \R$, we consider the differential
operators of first order (modified twistor operators)
\begin{displaymath}
\calp^s , \calq^t : \Gamma (\Sigma) \to \Gamma (TM \otimes \Sigma)
\end{displaymath}

defined by $\calp^s \psi := X^k \otimes \calp^s_{X_k} 
\psi , \calq^t \psi := X^k \otimes \calq^t_{X_k} \psi$ and
\begin{eqnarray*}
&& \calp^s_{X} \psi  :=  \cald_X \psi - s (\delta C(X) + \frac{1}{4n}
X \cdot dS) \cdot \psi , \\
&& \calq^t_{X}   :=  \cald_X \psi + t (\Ric - \frac{S}{n} )(X) \cdot
D \psi . 
\end{eqnarray*}

The image of $\cald$ is contained in the kernel of the Clifford multiplication,
i.e., 
\begin{equation} \label{gl-34}
X^k \cdot \cald_{X_k} \psi =0
\end{equation}

for all $\psi \in \Gamma (\Sigma)$. Thus, by (\ref{gl-10}) and 
(\ref{gl-25}), we see that the images of $\calp^s$ and $\calq^t$ are also 
contained in the kernel of the Clifford multiplication
\begin{equation} \label{gl-35}
X^k \cdot \calp^s_{X_k} \psi =0 \quad , \quad X^k \cdot
\calq^t_{X_k} \psi =0 . 
\end{equation}

For any $\psi \in \Gamma ( \Sigma)$, one has the well-known formula
\begin{equation} \label{gl-36}
| \cald \psi |^2 = | \nabla \psi |^2 - \frac{1}{n} | D \psi |^2 . 
\end{equation}

We introduce the selfadjoint nonnegative endomorphism
\begin{displaymath}
{\mathcal{E}} := E- \frac{1}{16n} |dS|^2 \stackrel{(\ref{gl-28})}{=}
F+ \frac{1}{16n (n-1)} |dS|^2
\end{displaymath}

and by straightforward calculations we obtain 
\begin{equation} \label{gl-37}
%\begin{array}{c}
| \calp^s \psi |^2 = | \cald \psi |^2 + 2s \mathrm{Re} \langle \delta
C(X^k) \nabla_{X_k} \psi , \psi \rangle +
 \frac{s}{2n} \mathrm{Re} \langle D \psi , dS \cdot \psi \rangle
+s^2 \langle \mathcal{E} \psi , \psi \rangle , 
%\end{array}
\end{equation}
\begin{equation} \label{gl-38}
%\begin{array}{c}
| \calq^t \psi |^2 = | \cald \psi |^2 - 2t \mathrm{Re} \langle \Ric (X^k)
\nabla_{X_k} \psi , D \psi \rangle +
 2t \frac{S}{n} | D \psi |^2 + t^2 | \Ric - \frac{S}{n} |^2 |
D  \psi |^2 . 
%\end{array}
\end{equation}

\vspace{0.3cm}

\begin{lem} \label{lem-3-1}
Let $\lambda$ be any eigenvalue of the Dirac operator $D$. Then, for 
all corresponding eigenspinors $\psi$ $(D \psi = \lambda \psi)$, it
holds that
\begin{eqnarray} \label{gl-39}
&& \frac{1}{2} (| \calp^t \psi |^2 + \calq^t \psi |^2 )= | \cald
\psi |^2 + t \frac{S}{n} \lambda^2 | \psi |^2 - \nonumber \\
&& - t(( \lambda^2 - \frac{S}{4} )( | \nabla \psi |^2 - (\lambda^2 - 
\frac{S}{4} | \psi |^2 ) + \frac{1}{4} | \Ric |^2 | \psi |^2 + \langle
\nabla_{\Ric (X_k)} \psi , \nabla_{X^k} \psi \rangle ) +\\
&&+ t \mathrm{div} (X_{\psi}) + \frac{t^2}{2} (\langle \mathcal{E} \psi , 
\psi \rangle + \lambda^2 | \Ric - \frac{S}{n}  |^2 | \psi |^2) , \nonumber
\end{eqnarray}

where $X_{\psi}$ is the vector field locally defined by 
\begin{displaymath}
X_{\psi} := \mathrm{Re} ( \langle ( D^2 - \frac{S}{4}) \psi , 
\nabla_{X^k} \psi \rangle + \langle \nabla_{X_j} D \psi + \frac{1}{2} \Ric
(X_j) \cdot \psi , X^k \cdot \nabla_{X^j} \psi \rangle ) X_k . 
\end{displaymath}
\end{lem}

\vspace{0.3cm}

\begin{proof}
By Lemma 1.4 in \cite{2} and (\ref{gl-24}), for all $\psi \in \Gamma (S)$, 
we have the identity
\begin{eqnarray} \label{gl-40}
&& \mathrm{Re} \langle \Ric (X^k) \nabla_{X_k} D \psi , \psi \rangle - 
\mathrm{Re} \langle \delta C(X^k) \nabla_{X_k} \psi , \psi \rangle =
\nonumber \\
&& | \nabla D \psi |^2 - | (D^2 - \frac{S}{4} ) \psi |^2 - \frac{S}{4}
| \nabla \psi |^2 + \frac{1}{4} | \Ric |^2 | \psi |^2 +\\
&& \langle \nabla_{\Ric (X^k)} \psi , \nabla_{X^k} \psi \rangle - 
\mathrm{div} (X_{\psi} ) . \nonumber
\end{eqnarray}

Using (\ref{gl-37}), (\ref{gl-38}) and (\ref{gl-40}) we obtain 
(\ref{gl-39}).
\end{proof}

\vspace{0.3cm}

Now, for $M$ being compact, 
let $\vartheta$ denote the supremum of all eigenvalues
of $\mathcal{E}$ on $\Sigma$. Then $\vartheta \ge 0$ and 
\begin{equation} \label{gl-42}
\langle \mathcal{E} \psi , \psi \rangle \le \vartheta | \psi |^2
\end{equation}

for any $\psi \in \Gamma (\Sigma)$. Moreover, let $\kappa_0$ be the infimum
of all eigenvalues of $\Ric$ on $TM$ and let $\kappa$ denote the 
supremum of its eigenvalues. Then, for any $\psi \in \Gamma (\Sigma)$, the
inequalities
\begin{equation} \label{gl-43}
\kappa_0 | \nabla \psi |^2 \le \langle \nabla_{\Ric (X_k)} \psi , 
\nabla_{X^k} \psi \rangle \le \kappa | \nabla \psi |^2
\end{equation}

are valid. We denote by $S_0$ the minimum of the scalar curvature $S$ and 
by $S_1$ its maximum and we use the notation
\begin{displaymath}
S_* := \left\{ \begin{array}{l}
S_0 \ \mbox{if $\kappa_0 \le 0$}\\[0.5em]
S_1 \ \mbox{if $\kappa_0 >0$ . }\end{array}  \right.
\end{displaymath}

Further, we introduce the functions $\alpha, \beta : \R \to \R$ defined by
\begin{eqnarray*}
\alpha (t) &:=& 1 + \frac{nt}{n-1} (\frac{S_1}{n} - \kappa_0 +
\frac{S_1 - S_0}{4}) + \frac{nt^2}{2(n-1)} | \Ric - \frac{S}{n} |^2_1 , \\
\beta (t) &:=& S_0 +t (| \Ric |^2_0 - S_* \kappa_0 + \frac{S_0 (S_1 - S_0)}
{4})- 2 \vartheta t^2 , 
\end{eqnarray*}

where $| \Ric |_0$ denotes the minimum of the function $| \Ric |$ and $|
\Ric - \frac{S}{n} |_1$ the maximum of $| \Ric - \frac{S}{n} |$.\\

\begin{theorem} \label{thm-3-2}
{\it 
Let $M$ be a compact Riemannian spin $n$-manifold and let $\lambda$ be any 
eigenvalue of the Dirac operator $D$. Then, for all $t \ge 0$, we have
\begin{equation} \label{gl-44}
\lambda^2 \ge \frac{n}{4(n-1)} \cdot \frac{\beta (t)}{\alpha (t)} . 
\end{equation}}
\end{theorem}

\vspace{0.3cm}

\begin{proof}
By Lemma 2.2 in \cite{6}, the inequalities
\begin{equation} \label{gl-45}
- \frac{S_1 - S_0}{4} (\lambda^2 - \frac{S_0}{4}) \int_M | \psi |^2 \le
\int_M ( \lambda^2 - \frac{S}{4})(| \nabla \psi |^2 - (\lambda^2 - \frac{S}{4})
| \psi |^2 \le \frac{S_1 - S_0}{4} (\lambda^2 - \frac{S_0}{4}) \int_M
| \psi |^2
\end{equation}

are valid for any eigenspinor $\psi$ to the eigenvalue
$\lambda$ of $D$. Using (\ref{gl-31}), (\ref{gl-36}), (\ref{gl-42}), and 
(\ref{gl-45}) we obtain (\ref{gl-44}) if we integrate the equation 
(\ref{gl-39}).
\end{proof}

\vspace{0.3cm}

We obtain the following corollary by computing the maximum of $\beta (t)$ for
$t \ge 0$.\\

\begin{cor} \label{cor-3-3}
There are no harmonic spinors on a compact Riemannian spin manifold with
$S_0 \le 0$ if the condition
\begin{equation} \label{gl-46}
| \Ric |^2_0 > S_0 ( \kappa_0 - \frac{S_1 -S_0}{4}) + \sqrt{8 | S_0 | 
\vartheta}
\end{equation}

is satisfied. In particular, the kernel of $D$ is trivial if $S_0 =0$ and 
$| \Ric |_0 >0$.
\end{cor}

\vspace{0.3cm}

\begin{NB} \label{NB-3-4}
(i) Our Corollary \ref{cor-3-3} is a generalization of Theorem 2.1 in 
\cite{2} since, in the case of a harmonic curvature tensor $(\delta R=0)$, we 
have $dS=0$ and $\vartheta =0$.\\
(ii) The inequality (\ref{gl-44}) can be written in the form
\begin{equation} \label{gl-47}
\lambda^2 \ge \frac{n}{4(n-1)} (S_0 +t \frac{\gamma (t)}{\alpha (t)}) , 
\end{equation}

where $\gamma (t)$ is the function given by
\begin{displaymath}
\gamma (t) := | \Ric |^2_0 - \frac{S_0}{n-1} (S_1 - \kappa_0 +
\frac{S_1 - S_0}{4}) - \kappa_0 (S_* - S_0) - 2t (\frac{nS_0}{4(n-1)}
| \Ric - \frac{S}{n} |^2_1 + \vartheta ) . 
\end{displaymath}

Thus, for $S_0 >0$, (\ref{gl-47}) yields a better estimate than 
(\ref{gl-01}) if $\gamma (t) >0$ for some $t>0$. We see immediately
that this is the case if the condition
\begin{equation} \label{gl-48}
| \Ric |^2_0 > \frac{S_0}{n-1} (S_1 - \kappa_0 + \frac{S_1 - S_0}{4} )
+ \kappa_0 (S_* - S_0)
\end{equation}

is fulfilled. This generalizes a corresponding assertion in \cite{2}, 
Section 2. In particular, if $S$ is constant and positive, (\ref{gl-48}) 
simplifies to
\begin{equation} \label{gl-49}
| \Ric |^2_0 > \frac{S}{n-1} (S - \kappa_0 ) . 
\end{equation}

(iii) The limiting case of (\ref{gl-44}) corresponds to the limiting case of 
(\ref{gl-01}) since, by the same arguments that we used in Section 2 of 
\cite{2}, it follows that (\ref{gl-44}) can be an equality for the first
eigenvalue of $D$ for $t=0$ only.
\end{NB}

In order to write down the main result of this section the notations
\begin{eqnarray*}
&& A :=  | \Ric |^2_0 - \frac{S_0}{n-1} (S_1 - \kappa_0 + \frac{S_1 - S_0}{4} )
- \kappa_0 (S_* - S_0) , \\
&& b := \frac{n}{n-1} ( \frac{S_1}{n} - \kappa_0 + \frac{S_1 - S_0}{4}) 
\quad , \quad c := | \Ric - \frac{S}{n} |_1 \sqrt{\frac{2n}{n-1}} , \\
&& a := \frac{4}{A} ( \frac{n S_0}{4(n-1)}  | \Ric - \frac{S}{n} |^2_1 + 
\vartheta )
\end{eqnarray*}

are convenient. The function $t \gamma (t) / \alpha (t)$ attains its
maximum for $t >0$ if the condition (\ref{gl-48}) is satisfied, i.e., 
if $A>0$. By computing this maximum and assertion (iii) of Remark
\ref{NB-3-4}, we obtain the following result.\\

\begin{cor} \label{thm-3-5}
{
Let $M$ be a compact Riemannian spin $n$-manifold with $A \ge 0$. Then, for 
every eigenvalue $\lambda$ of the Dirac operator, we have the inequality
\begin{equation} \label{gl-50}
\lambda^2 \ge \frac{n}{4(n-1)} (S_0 + \frac{A}{a+b+ \sqrt{a^2+2ab +c^2}}) , 
\end{equation}

which is never an equality if $A>0$.}
\end{cor}

\vspace{0.3cm}

\begin{cor} \label{cor-3-6}
If $M$ is a compact Riemannian spin $n$-manifold such that $S_0 =0$ and
$| \Ric |_0 >0$, then every eigenvalue $\lambda$ of the Dirac operator
satisfies the estimate 
\begin{equation} \label{gl-51}
\lambda^2 > \frac{n}{4(n-1)} \cdot \frac{| \Ric |^2_0}{a+b+
\sqrt{a^2 +2ab +c^2}}
\end{equation}

with the constants
\begin{displaymath}
a= \frac{4 \vartheta}{| \Ric |^2_0} \quad , \quad b= \frac{n}{n-1}
(\frac{n+4}{4n} S_1 - \kappa_0) \quad , \quad c= | \Ric - 
\frac{S}{n} |_1 \sqrt{\frac{2n}{n-1}} . 
\end{displaymath}
\end{cor}

\vspace{0.3cm}
 
\begin{NB} \label{NB-3-7}
(i) Our Corollary \ref{thm-3-5} is comparable with Theorem 3.1. in \cite{2},
which uses the additional assumption that 
$\delta R=0$. But Corollary \ref{thm-3-5}
is not a direct generalization of this Theorem 3.1. since the application
of Corollary \ref{thm-3-5} to the case of a harmonic curvature tensor yields
a weaker result than Theorem 3.1. In particular, applying Corollary
\ref{cor-3-6} to the special case of $\delta R=0$, the estimate (\ref{gl-51})
may be written as
\begin{equation} \label{gl-52}
\lambda^2 > \frac{1}{4} \cdot \frac{| \Ric |^2_0}{| \kappa_0 | + | \Ric |_1
\sqrt{\frac{2(n-1)}{n}}}
\end{equation}

since $\delta R=0$ implies $E=0$ and $dS=0$ and, hence, $\vartheta =0$. 
Comparing (\ref{gl-03}) and (\ref{gl-52}) we see that (\ref{gl-52}) is a 
weaker estimate than (\ref{gl-03}).\\
(ii) Corollary 4.1 in \cite{6} is a result similar to Corollary 
\ref{cor-3-6}, it  was obtained under the additional assumption 
that the Ricci 
tensor commutes with its covariant derivatives of first order $([\Ric, 
\nabla_X \Ric ]=0)$.\\
(iii) The Examples 4.1. and 4.2. in \cite{6} yield simple examples of
manifolds for which the lower bounds in the estimates (\ref{gl-50}) or
(\ref{gl-51}), respectively, can be computed easily.
\end{NB}

\vspace{0.3cm}

\section{\bf Weyl tensor depending estimates}
 
Our estimate (\ref{gl-50}) cannot be better than 
(\ref{gl-01}) if $M$ is Einstein or if $| \Ric |_0
=0$. In this section we prove estimates that also
work in such situations. For $s,t \in \R$, let
\begin{displaymath}
\calr^s , \cals^t : \Gamma (\Sigma) \to
\Gamma (TM \otimes \Sigma )
\end{displaymath}

be the first order differential operators locally
defined by $\calr^s \psi := X^k \otimes \calr^s_{X_k}
\psi , \cals^t \psi = X^k \otimes \cals^t_{X_k} \psi$ with
\begin{displaymath}
\calr^s_X \psi := \cald_X \psi - s \delta B(X) \psi \quad ,
\quad \cals^t_X \psi := \cald_X \psi - t B(X, X^k)
\nabla_{X_k} \psi . 
\end{displaymath}

Then, for any $\psi \in \Gamma (\Sigma)$, we have
\begin{equation} \label{gl-53}
| \calr^s \psi |^2 = | \cald \psi |^2 +2s
\mathrm{Re} \langle \delta B(X^k) \nabla_{X_k} \psi,
\psi \rangle + s^2 \langle F \psi , \psi \rangle , 
\end{equation}
\begin{equation} \label{gl-54}
\begin{array}{l}
| \cals^t \psi |^2 = | \cald \psi |^2 - 2t
\mathrm{Re} \langle \delta B(X^k) \nabla_{X_k} \psi , 
\psi \rangle - t \langle H \psi , \psi \rangle +\\[0.8em]
+2t \mathrm{div} (\mathrm{Re} \langle B(X^k , X^l)
\nabla_{X_l} \psi , \psi \rangle X_k) +t^2 \langle
B^2 (X^k , X^l) \nabla_{X_k} \psi , \nabla_{X_l}
\psi \rangle 
\end{array}
\end{equation}

and, hence,
\begin{eqnarray} \label{gl-55}
&& \frac{1}{2} (| \calr^{2t} \psi |^2 + | \cals^{2t}
\psi |^2 )= | \cald \psi |^2 - t \langle H \psi , \psi \rangle
+ \nonumber \\[0.7em]
&& + 2t \mathrm{div} (\mathrm{Re} \langle B(X^k, X^l)
\nabla_{X_l} \psi , \psi \rangle X_k )+\\[0.7em]
&& +2t^2 ( \langle  F \psi , \psi \rangle + \langle
B^2 (X^k , X^l) \nabla_{X_k} \psi , \nabla_{X_l}
\psi \rangle ) . \nonumber
\end{eqnarray}

\vspace{0.3cm}

\begin{theorem} \label{thm-4-1}
{\it Let $M$ be a compact Riemannian spin $n$-manifold
with harmonic Weyl tensor $(\delta W=0)$ and let
$\lambda$ be any eigenvalue of the Dirac operator.
Then, for all $t \ge 0$, the inequality
\begin{equation} \label{gl-56}
\lambda^2 \ge \frac{n}{4(n-1)} (S_0 + 
\frac{4 \nu_0 t - (n-1) \mu^2 S_0 t^2}{1+ n(n-1)
\mu^2 t^2} )
\end{equation}

is valid.}
\end{theorem}

\vspace{0.3cm}

\begin{proof}
By (\ref{gl-21}), $\delta W=0$ implies $\delta B=0$. Integrating
equation (\ref{gl-54}) for any eigenspinor $\psi (D \psi =
\lambda \psi)$ we find (\ref{gl-56}) by using $\delta B=0$, 
(\ref{gl-30}), (\ref{gl-31}), (\ref{gl-33}) and (\ref{gl-36}).
\end{proof}

\vspace{0.3cm}

The following result is proved by computing
the maximum of the right-hand side of (\ref{gl-56}) for $t \ge 0$.\\

\begin{cor} \label{thm-4-2}
{\it Let $M$ be a compact Riemannian spin $n$-manifold with
$\delta W =0$ and $\mu >0$. Then every eigenvalue $\lambda$
of the Dirac operator satisfies the estimate
\begin{equation} \label{gl-57}
\lambda^2 \ge \frac{1}{8(n-1)} \Big((2n-1) S_0 +
\sqrt{S^2_0 + \frac{n}{n-1} (\frac{4 \nu_0}{\mu} )^2}\Big) . 
\end{equation}

For $S_0 \le 0$, this lower bound is positive if
\begin{equation} \label{gl-58}
\nu_0 > \frac{n-1}{2} | S_0 | \mu . 
\end{equation}

In particular, there are no harmonic spinors if $S_0 =0$
and $\nu_0 >0$.}
\end{cor}

\vspace{0.3cm}

Every Einstein manifold fulfils the condition $\delta
W=0$. Thus, we obtain\\

\begin{cor} \label{cor-4-3}
The estimate (\ref{gl-57}) is valid on any compact
Einstein spin manifold with $\mu >0$.
\end{cor}

\vspace{0.3cm}

\begin{NB} \label{NB-4-4}
(i) Comparing (\ref{gl-07}) and (\ref{gl-57}) we see that
(\ref{gl-57}) is the better estimate. Thus, our
Corollary 4.1 improves Theorem 3.1 in \cite{3}.\\
(ii) For $S_0 >0$, (\ref{gl-57}) yields a better estimate
than (\ref{gl-01}) if $\nu_0 >0$. By Corollary 4.2,
this is also the case if the manifold is Einstein or even 
Ricci flat.
\end{NB}

\vspace{0.3cm}

Our next aim is to prove an estimate similar to (\ref{gl-57})
for manifolds whose Weyl tensor is not harmonic. We denote by 
$\eta$ the supremum of all eigenvalues of the endomorphism
$F$ on $\Sigma$. Then $\eta \ge 0$ and it holds that
\begin{equation} \label{gl-59}
\langle F \psi , \psi \rangle \le \eta | \psi |^2
\end{equation}

for all $\psi \in \Gamma (\Sigma)$.\\

\begin{theorem} \label{thm-4-5}
{\it Let $M$ be any compact Riemannian spin $n$-manifold
and let $\lambda$ be any eigenvalue of the Dirac operator.
Then, for all $t \ge 0$, we have the inequality
\begin{equation}  \label{gl-60}
\lambda^2 \ge \frac{n}{4(n-1)} (S_0 + \frac{4 \nu_0 t-
2((n-1) \mu^2 S_0 + 4 \eta) t^2}{1+2n (n-1) \mu^2 t^2}
) . 
\end{equation}}
\end{theorem}

\vspace{0.3cm}

\begin{proof}
Using (\ref{gl-59}) we integrate the equation 
(\ref{gl-55}) and find (\ref{gl-60}) by simple estimates
as before.
\end{proof}

\vspace{0.3cm}

By computing the maximum of the right-hand side of (\ref{gl-60})
with respect to $t \ge 0$, we obtain the following  result.\\

\begin{cor} \label{thm-4-6}
{If $M$ is a compact Riemannian spin $n$-manifold with 
$\mu >0$, then, for every eigenvalue $\lambda$ of the Dirac
operator $D$, the estimate
\begin{equation} \label{gl-61}
\lambda^2 \ge \frac{1}{8(n-1)}\Big((2n-1) S_0 - 
\frac{4 \eta}{(n-1) \mu^2} + \sqrt{(S_0 + \frac{4 \eta}{(n-
1) \mu^2} )^2 + \frac{8n}{n-1} (\frac{\nu_0}{\mu})^2} \Big)
\end{equation}

is valid. For $S_0 \le 0$, this lower bound is positive
and, hence, $\ker (D)=0$ if the condition
\begin{equation} \label{gl-62}
\nu_0 > \sqrt{| S_0 | (2 \eta + \frac{1}{2} (n-1)^2
\mu^2 | S_0 | )} 
\end{equation}

is fulfilled.}
\end{cor}

\vspace{0.3cm}

\begin{cor} \label{cor-4-7}
For every eigenvalue $\lambda$ of the Dirac operator
on a compact Riemannian spin $n$-manifold with $S_0 =0$
and $\nu_0 >0$, we have the estimate
\begin{equation} \label{gl-63}
\lambda^2 \ge \frac{n}{4(n-1)} \cdot \frac{\nu_0^2}{\eta
+ \sqrt{\eta^2 + (\begin{array}{c} n\\2 \end{array} )
\mu^2 \nu^2_0}} . 
\end{equation}

In particular, there are no harmonic spinors.
\end{cor}

\vspace{0.3cm}

\begin{NB} \label{NB-4-8}
                     
(i) For $S_0 >0$, (\ref{gl-61}) also yields a better
estimate than (\ref{gl-01}) if $\nu_0 >0$.\\
(ii) It is not known if there exist manifolds with the 
property that (\ref{gl-57}) or (\ref{gl-61}), 
respectively, is an equality for the first eigenvalue
$\lambda_1$ of the Dirac operator.
\end{NB}

\vspace{0.3cm}

\section{\bf Estimates depending on the whole curvature
tensor}

In order to obtain estimates for the first eigenvalue of
the Dirac operator that depend on the Ricci tensor and
also on the Weyl tensor we consider, for all $t \in \R$,
the first order differential operator
\begin{displaymath}
\calt^t : \Gamma (\Sigma) \to \Gamma (TM \otimes \Sigma) , 
\end{displaymath}

which is locally defined by $\calt^t \psi := X^k \otimes \calt^t_{X_k}
\psi$, and 
\begin{displaymath}
\calt^t_X \psi := \cald_X \psi - tC (X, X^k) \nabla_{X_k} \psi.
\end{displaymath}

Then, for any $\psi \in \Gamma (\Sigma)$, it holds that
\begin{eqnarray} \label{gl-64}
&& | \calt^t \psi |^2 = | \cald \psi |^2 + \frac{t}{n}
\mathrm{Re} \langle \Ric (X^k) \nabla_{X_k} \psi , D \psi
\rangle 
 -2t \mathrm{Re} \langle \delta C(X^k) \nabla_{X_k} \psi , \psi
\rangle -t \langle G \psi , \psi \rangle \\
&& -2t \mathrm{div} (\mathrm{Re} \langle C (X^k , X^l)
\nabla_{X_l} \psi , \psi \rangle X_k ) + t^2 \langle C^2 (
X^k , X^l ) \nabla_{X_k} \psi , \nabla_{X_l} \psi \rangle . \nonumber
\end{eqnarray}

\vspace{0.3cm}

\begin{lem} \label{lem-5-1}
Let $M$ be a Riemannian spin $n$-manifold and let $\lambda$ be any
eigenvalue of the Dirac operator $D$. Then, for any corresponding 
eigenspinor $\psi (D \psi = \lambda \psi)$ and all $t \in \R$, we have
the equations
\begin{eqnarray} \label{gl-65} 
&& | \calt^t \psi |^2 = | \cald \psi |^2 - t \frac{2n-1}{n}
\mathrm{Re} \langle \delta C(X^k) \nabla_{X_k} \psi , \psi
\rangle \nonumber \\
&& + \frac{t}{n} (( \lambda^2 - \frac{S}{4} )( | \nabla \psi |^2 - 
(\lambda^2 - \frac{S}{4} ) | \psi |^2 )+ \langle \nabla_{\Ric (X_k)}
\psi , \nabla_{X^k} \psi \rangle ) \nonumber \\
&& -t ( \langle H \psi , \psi \rangle + \frac{1}{4n} ( \frac{n+2}{n
-2} | \Ric - \frac{S}{n} |^2 + \frac{S^2}{n(n-1)} ) | \psi |^2)\\
&& \mathrm{div} ( \frac{t}{n} X_{\psi} + 2t \mathrm{Re} \langle 
C(X^k , X^l) \nabla_{X_l} \psi , \psi \rangle X_k )
 + t^2 \langle C^2 (X^k , X^l ) \nabla_{X_k} \psi , \nabla_{X_l}
\psi \rangle , \nonumber \\[1em]
&& \frac{1}{2} (| \calp^{\frac{2n-1}{n}t} \psi |^2 + | \calt^{2t} \psi 
|^2 )= | \cald \psi |^2 \nonumber \\
&& + \frac{t}{n} (( \lambda^2 - \frac{S}{4})( | \nabla \psi |^2 - 
(\lambda^2 - \frac{S}{4}) | \psi |^2 ) + \langle \nabla_{\Ric (X_k)} \psi , 
\nabla_{X^k} \psi \rangle ) \nonumber\\
&& 
-t ( \langle H \psi , \psi \rangle + \frac{1}{4n} ( \frac{n+2}{n-2} | \Ric
- \frac{S}{n} |^2 + \frac{S^2}{n(n-1)} ) | \psi |^2  \label{gl-66} \\
&& - \mathrm{div} ( \frac{t}{n} X_{\psi} + 2t \mathrm{Re} \langle C(
X^k , X^l) \nabla_{X_l} \psi , \psi \rangle X_k ) \nonumber \\
&& +2t^2 ( \langle C^2 (X^k , X^l) \nabla_{X_k} \psi , \nabla_{X_l}
\psi \rangle + ( \frac{2n-1}{2n} )^2 \langle \cale \psi , \psi \rangle ).
\nonumber
\end{eqnarray}
\end{lem}

\vspace{0.3cm}

\begin{proof}
Inserting (\ref{gl-19}) and (\ref{gl-40}) into (\ref{gl-64}) we find
(\ref{gl-65}). Using (\ref{gl-37}) and (\ref{gl-65}) we obtain
(65).
\end{proof}

\vspace{0.3cm}

Again, let $M$ be compact. By $| S |_0$ we denote the minimum of the
function $|S|$ on $M$ and we use the notation
\begin{displaymath}
S_{\star} := \left\{ \begin{array}{ll}
S_0 & \mbox{if $\kappa \ge 0$}\\
S_1 & \mbox{if $\kappa <0$} \end{array} \right. . 
\end{displaymath}

Moreover, we introduce six functions $\alpha_p , \beta_p , \gamma_p :
\R \to \R , p \in \{ 1,2 \}$, defined by
\begin{eqnarray*}
&&\alpha_p (t) := 1 + \frac{t}{n-1} ( \kappa + \frac{S_1 - S_0}{4} ) +pn
(n-1) \zeta^2 t^2 , \\
&&\beta_p (t) := S_0 + t( 4 \nu_0 + \frac{1}{n} ( \frac{n+2}{n-2} |
\Ric -  \frac{S}{n} |^2_0 + \frac{| S|^2_0}{n(n-1)} + \frac{S_0 (S_1 - 
S_0)}{4} + S_{\star} \kappa ))\\
&& \mbox{} \hspace{1.2cm}
+pt^2 ((n-1)^2 S_0 \zeta^2 - ( \frac{2n-1}{n})^2 \vartheta ) , \\
&& \gamma_p (t) := 4 \nu_0 + \frac{1}{n} (\frac{n+2}{n-2} | \Ric - 
\frac{S}{n} |^2_0 + \frac{| S|^2_0}{n(n-1)} - \frac{S_0}{n-1} ( \kappa
+\frac{S_1 - S_0}{4} ) + \kappa (S_{\star} - S_0 ))\\
&& \mbox{} \hspace{1.2cm}
-pt (( n-1) S_0 \zeta^2 + ( \frac{2n-1}{n} )^2 \vartheta) .
\end{eqnarray*}

\vspace{0.3cm}
 
\begin{theorem} \label{thm-5-2}
{\it Let $\lambda$ be any eigenvalue of the Dirac operator on a compact
Riemannian spin $n$-manifold. Then the following holds: \\
(i) For any $t \ge 0$ with $\beta_2 (t) >0$, we have the estimate
\begin{equation} \label{gl-67}
\lambda^2 \ge \frac{n}{4(n-1)} \cdot \frac{\beta_2 (t)}{\alpha_2
(t)} = \frac{n}{4(n-1)} (S_0 + t \frac{\gamma_2 (t)}{\alpha_2 (t)}) . 
\end{equation}

(ii) If the curvature tensor is harmonic, then the estimate
\begin{equation} \label{gl-68}
\lambda^2 \ge \frac{n}{4(n-1)} \cdot \frac{\beta_1 (t)}{\alpha_1 (t)}
= \frac{n}{4(n-1)} (S +t \frac{\gamma_1 (t)}{\alpha_1 (t)} )
\end{equation}

is valid for every $t \ge 0$ with $\beta_1 (t) >0$. }
\end{theorem}

\vspace{0.3cm}

\begin{proof}
Integrating equation (\ref{gl-66}) and using (\ref{gl-33}), (\ref{gl-42})
and (\ref{gl-43}), for any $t \ge 0$, we obtain
\begin{equation} \label{gl-69}
\lambda^2 \alpha_2 (t) \ge \frac{n}{4(n-1)} \beta_2 (t) = \frac{n}{4(n-1)}
(S_0 \alpha_2 (t) + \gamma_2 (t)) . 
\end{equation}

In particular, (\ref{gl-69}) shows that $\beta_2 (t) >0 \ (t \ge 0)$ 
forces $\alpha_2 (t) >0$. 
This proves the assertion (i) of our theorem.
Further, the supposition $\delta R=0$ implies $\delta C=0$ by 
(\ref{gl-21}) and, moreover, $\vartheta =0, S_0 =S_1 =S$. Thus, integrating
equation (\ref{gl-65}) one analogously proves the assertion (ii).
\end{proof}

\vspace{0.3cm}

\begin{cor} \label{cor-5-3}
On a compact Riemannian spin $n$-manifold with $S_0 \le 0$, we have the 
following:\\
(i) There are no harmonic spinors if the condition
\begin{equation} \label{gl-70}
\begin{array}{l}
\displaystyle 4 n \nu_0 + \frac{n+2}{n-2} | \Ric - \frac{S}{n} |^2_0 +
\frac{| S|^2_0}{n(n-1)} + S_{\star} \kappa >\\[0.7em]
\displaystyle 
\frac{| S_0| (S_1 - S_0)}{4} + 4 \sqrt{2 | S_0| (( \begin{array}{c}
n\\2 						   \end{array} )^2
| S_0 | \zeta^2 + (\frac{2n-1}{2})^2 \vartheta )}
\end{array}
\end{equation}

is satisfied. In particular, for $S_0 =0$, there are no harmonic spinors
if $\nu_0 >0$ or $|\Ric - \frac{S}{n} |_0 >0$.\\
(ii) If the curvature tensor is harmonic, then there exist no harmonic
spinors if
\begin{equation} \label{gl-71}
4n \nu_0 + \frac{n+2}{n-2} | \Ric - \frac{S}{n} |^2_0 + \frac{S^2}{n(n-1)}
> |S| (\kappa + 4 ( \begin{array}{c} n\\2 \end{array} ) \zeta ) . 
\end{equation}

In particular, for $S=0$, we have $\ker (D) =0$ if $\nu_0 >0$ or 
$| \Ric |_0 >0$.
\end{cor}

\vspace{0.3cm}

\begin{proof}
(\ref{gl-70}) implies that the function $\beta_2 (t)$ attains positive
values for some $t >0$. The condition (\ref{gl-71}) implies that also 
the function $\beta_1 (t)$ has this property.
\end{proof}

\vspace{0.3cm}

\begin{NB} \label{NB-5-4}
(i) If the condition
\begin{equation} \label{gl-72}
4n \nu_0 + \frac{n+2}{n-2} | \Ric - \frac{S}{n} |^2_0 +
\frac{| S  |^2_0}{n(n-1)} > \frac{S_0}{n-1} (\kappa + \frac{S_1 - S_0}{4})
\end{equation}

is satisfied on a compact Riemannian spin $n$-manifold with $S_0 >0$, then
(\ref{gl-67}) yields a better estimate than (\ref{gl-01}) since this 
condition implies that the function $\gamma_2 (t)$ attains  positive
values for some $t >0$. We note that $S_0 >0$ implies $\kappa >0$ and, 
hence, $\alpha_2 (t) \ge 1$ for $t \ge 0$.\\
(ii) In the case of a harmonic curvature tensor, the function $\gamma_1 (t)$ 
reaches positive values for some $t >0$ if 
\begin{equation} \label{gl-73}
4n \nu_0 + \frac{n+2}{n-2} | \Ric - \frac{S}{n} |^2_0 > \frac{S}{n-1}
(\kappa - \frac{S}{n} ) . 
\end{equation}

Thus, if $S >0$ and (\ref{gl-73}) is fulfilled, (\ref{gl-68}) yields a
better estimate than (\ref{gl-01}).\\
(iii) The assertion (ii) of Corollary \ref{cor-5-3} is an improvement
of the Theorem 4.1 in \cite{3}, where, instead of $\zeta$, another curvature
invariant $\sigma$ was used. $\zeta$ and $\sigma$ are related by 
\begin{equation} \label{gl-74}
\zeta \le \frac{1}{2} (\begin{array}{c} n\\2 \end{array} ) \sigma
\end{equation}

(see \cite{3}, Section 4). Replacing $\zeta$ by the value 
$(\begin{array}{c} n \\2 \end{array} )\sigma /2$ inequality
(\ref{gl-71}) becomes
a condition that is weaker than the condition (38) in \cite{3}.
\end{NB}

\vspace{0.3cm}

In the end of this paper we show that another combination of our basic
Weitzenb\"ock formulas leads to similar results, but 
they do not contain
the curvature invariants $\kappa_0$ and $\kappa$. Using (\ref{gl-37}),
(\ref{gl-38}) and (\ref{gl-64}) we find the equation
\begin{eqnarray} \label{gl-75}
&& \frac{1}{2} (| \calq^{\frac{t}{n}} \psi |^2 + | \calt^{2t} \psi |^2 )=
| \cald \psi |^2 + t \frac{S}{n} | D \psi |^2 \nonumber \\
&& -t \langle G \psi , \psi \rangle -2t \mathrm{Re} \langle \delta C(
X^k) \nabla_{X_k} \psi , \psi \rangle - 2t \mathrm{div} (\mathrm{Re}
\langle C(X^k , X^l) \nabla_{X_l} \psi , \psi \rangle X_k )\\
&& +2t^2 (\frac{1}{4n^2} | \Ric - \frac{S}{n} |^2 | D \psi |^2 + 
\langle C^2 (X^k , X^l) \nabla_{X_k} \psi , \nabla_{X_l} \psi \rangle )
\nonumber
\end{eqnarray}

and, moreover, 
\begin{equation} \label{gl-76}
\begin{array}{l} \displaystyle
\frac{1}{3} (| \calp^{3t} \psi |^2 + | \calq^{\frac{3t}{2n}} \psi |^2
+ | \calt^{3t} \psi |^2 )=\\[0.5em]
\displaystyle
| \cald \psi |^2 - t \langle G \psi , \psi \rangle  +t \frac{S}{n^2}
| D \psi |^2 \\[0.5em]
\displaystyle 
+ \frac{t}{2n} \mathrm{Re} \langle  D \psi , dS \cdot \psi \rangle -
2t \mathrm{div} (\mathrm{Re} \langle C(X^k , X^l ) \nabla_{X_l} \psi , 
\psi \rangle X_k ) \\[0.5em]
\displaystyle 
+3t^2 ( \langle \cale \psi , \psi \rangle + \frac{1}{4n^2} | \Ric -
\frac{S}{n} |^2 | D \psi |^2 + \langle C^2 (X^k , X^l ) \nabla_{X_k}
\psi , \nabla_{X_l} \psi \rangle ) . 
\end{array}
\end{equation}

Both equations are valid for any $t \in \R$ and any $\psi \in \Gamma 
(\Sigma)$. We introduce the six functions $\alpha_p, \beta_p, \gamma_p :
\R \to \R , p \in \{ 3,4 \}$, defined by\\

$\displaystyle 
\alpha_p (t) := 1+t \frac{S_1}{n(n-1)} + (p-1) t^2 ( \frac{1}{4n(n-1)}
| \Ric - \frac{S}{n} |^2_1 + n (n-1) \zeta^2 ) ,$ \\

$\displaystyle
\beta_p (t) := S_0 +t ( 4 \nu_0 + \frac{2}{n-2} | \Ric - \frac{S}{n} |^2_0
+ \frac{|S|^2_0}{n(n-1)} ) +(p-1) t^2 ((n-1)^2 S_0 \zeta^2 - 4 \vartheta),$\\

$\displaystyle 
\gamma_p (t) := 4 \nu_0 + \frac{2}{n-2} | \Ric - \frac{S}{n} |^2_0 -
\frac{S_0}{n(n-1)} (S_1 - S_0) - (p-1) tS_0 ( \frac{1}{4n(n-1)} | \Ric
- \frac{S}{n} |^2_1 + (n-1) \zeta^2 + 4 \vartheta ).  $\\

\vspace{0.3cm}

\begin{theorem} \label{thm-5-5}
{\it Let $\lambda$ be any eigenvalue of the Dirac operator on a compact
Riemannian spin $n$-manifold. Then the following holds:\\
(i) For every $t \ge 0$ with $\beta_4 (t) >0$, we have the estimate
\begin{equation} \label{gl-77}
\lambda^2 \ge \frac{n}{4(n-1)} \cdot \frac{\beta_4 (t)}{\alpha_4 (t)} =
\frac{n}{4(n-1)} (S_0 + t \frac{\gamma_4 (t)}{\alpha_4 (t)}).
\end{equation}
(ii) In the special case that $\delta R=0$, the estimate
\begin{equation} \label{gl-78}
\lambda^2 \ge \frac{n}{4(n-1)} \cdot \frac{\beta_3 (t)}{\alpha_3 (t)}
= \frac{n}{4(n-1)} (S+ t \frac{\gamma_3 (t)}{\alpha_3 (t)})
\end{equation}

is valid for every $t \ge 0$ with $\beta_3 (t) >0$. }
\end{theorem}

\vspace{0.3cm}

\begin{proof}
Inserting any eigenspinor $\psi$ to the eigenvalue $\lambda$ of $D$
into equation (\ref{gl-76}) and then integrating this equation we obtain
(\ref{gl-77}) by (\ref{gl-19}), (\ref{gl-32}), (\ref{gl-42}) and 
analogous considerations as in the proof of Theorem \ref{thm-5-2}. In
the special case of $\delta R=0$, we integrate equation (\ref{gl-75})
for any eigenspinor $\psi$. Then we find (\ref{gl-78}).
\end{proof}

\vspace{0.3cm}

Studying the conditions under which the functions $\beta_3 (t)$ and 
$\beta_4 (t)$, respectively, attain positive values for some $t >0$, we
immediately obtain the next result.\\

\begin{cor} \label{cor-5-6}
The following holds on a compact Riemannian spin $n$-manifold with
$S_0 \le 0$:\\
(i) There are no harmonic spinors if
\begin{equation} \label{gl-79}
4n \nu_0 + \frac{2n}{n-2} | \Ric - \frac{S}{n} |^2_0 + \frac{| S|^2_0}{
n-1} > 4 \sqrt{ 3 | S_0 | (( \begin{array}{c} n\\2 \end{array} )^2 | S_0
| \zeta^2 + n^2 \vartheta)} . 
\end{equation}

In particular, for $S_0 =0$, we have $\ker (D) =0$ if $\nu_0 >0$ or
$| \Ric - \frac{S}{n} |_0 >0$.\\
(ii) In the special situation that $\delta R=0$, there are no harmonic spinors
if
\begin{equation} \label{gl-80}
4n \nu_0 + \frac{2n}{n-2} | \Ric - \frac{S}{n} |^2_0 + \frac{S^2}{n-1}
> 4 (\begin{array}{c} n\\2 \end{array}) \zeta |S| \sqrt{2} . 
\end{equation}
\end{cor}

\vspace{0.3cm}

\begin{NB} \label{NB-5-7}
(i) For $S_0 >0$, (\ref{gl-77}) gives a better estimate than (\ref{gl-01})
if
\begin{equation} \label{gl-81}
4n \nu_0 + \frac{2n}{n-2} | \Ric - \frac{S}{n} |^2_0 > \frac{S_0}{n-1}
(S_1 - S_0) . 
\end{equation}

(ii) In the special case of a harmonic curvature tensor and $S>0$, 
(\ref{gl-78}) yields a better estimate than (\ref{gl-01}) if $\nu_0 
>0$ or $| \Ric - \frac{S}{n} |_0 >0$.\\
(iii) The same arguments that are used in the proof of Theorem 4.2
in \cite{3} show that, for an optimal parameter $t_0 >0$, the inequalities
(\ref{gl-67}), (\ref{gl-68}) and (\ref{gl-77}), (\ref{gl-78}) can never
be equalities for the first eigenvalue of the Dirac operator.\\
(iv) If the first order covariant derivatives of the Ricci tensor commute
$([ \nabla_X \Ric , \nabla_Y \Ric ] =0)$, we see, by (\ref{gl-29}), that the 
number $\vartheta$, which enters the estimates (\ref{gl-44}), (\ref{gl-67})
and (\ref{gl-77}), is simply the maximum of the function $\frac{1}{4} |
\nabla \Ric |^2 - \frac{n+1}{16n} |dS |^2$. Moreover, 
in this case it becomes obvious, owing to (\ref{gl-28}) 
that the number $\eta$, which occurs in Section 4,
is given by the maximum of the function $\frac{1}{4} | \nabla \Ric |^2 - 
\frac{n}{16(n-1)} |dS |^2$ then.
\end{NB}

%\vspace{0.3cm}

\vspace{1cm}

%\refname{References}

\end{document}